\documentclass[11pt,draft]{amsart}
\usepackage{amssymb,latexsym}

\usepackage{setspace}

\theoremstyle{plain}
\newtheorem{theorem}{Theorem}

\newtheorem{definition}{Definition}

\theoremstyle{definition}

\newcommand{\curl}{\nabla \times } 
\newcommand{\dvg}{\nabla \cdot} 
\begin{document}

\onehalfspacing

\title[steady flow of nonhomogeneous asymmetric fluid]{On the steady
viscous flow of a nonhomogeneous asymmetric fluid}
\author{F\'abio Vitoriano e Silva}
\address{IME-UFG\\
         Caixa Postal 131\\
         74001-970 Goi\^ania, Goi\'as, Brazil}
\email{fabio@mat.ufg.br}
\thanks{Submitted to Annali di Mat. Pura ed Applicata on May 19, 2011}

\date{\today}

\begin{abstract}
We consider a boundary value problem for
the system of equations describing the stationary motion of a viscous
nonhomogeneous asymmetric fluid in a bounded planar domain
having a $C^2$ boundary. We use a stream-function
formulation after the manner of N. N. Frolov [Math. Notes, \textbf{53}(5-6), 650--656, 1993] in which the fluid density depends on the stream-function by means of another function determined by the boundary conditions. This allows for dropping some of the equations, most notably the continuity equation. With a fixed point argument we show the existence of solutions to the resulting system.
\end{abstract}
\maketitle

 \section*{Introduction}

Density dependent fluids are well known and have been studied by several authors.  Antontsev, Kazhikov and Monakhov treat in~\cite{akm93} an assortment of problems on density dependent flows of either compressible or incompressible Newtonian fluids. A more recent account on such problems and some improvements on results in~\cite{akm93} are available in~\cite{lions96} and references therein. 

The case of non-Newtonian fluids is less studied then the previous ones.  Several constitutive laws lead to such fluids and, among them, we are particularly interested in the non-symmetric fluids, named micropolar fluids, introduced in~\cite{er_66}. Particles of these fluids undergo translations and rotations as well and their theory has proved to be useful in describing phenomena in which the structure of the fluid should be accounted for, \emph{e.g.},~blood flow in thin vessels or flows of some slurries and polymeric fluids, see~\cite{lu_99}. 

Some authors studied evolutionary density dependent flows of micropolar fluids, such as the papers~\cite{bmkc03,cgmk02,guillenetal10} whereas, to the best of our knowledge, basic results on flows in a stationary regime are still lacking. Our aim in this paper is, therefore, to extend the result in~\cite{frolov93} placing the theory of micropolar fluids in a similar level of knowledge as the theory of the standard Navier-Stokes fluids. 

Frolov addressed in~\cite{frolov93} the stationary 2D flow of an incompres\-si\-ble inhomogeneous fluid resorting to a judicious stream-function formulation. 
His formulation has been successfully adapted to solve new problems as the mixing of two fluids having different (and discontinuous) densities~\cite{msantos} and to the boundary control of steady state flow of a viscous incompressible nonhomogeneous fluid~\cite{illa01}.

The field equations of the model, in a steady state regime, form the following system
\begin{equation}\label{e:vect_f}
 \begin{aligned}
- (\mu+\mu_r)\Delta \vec v +\rho (\vec v \cdot \nabla) \vec v &+ \nabla p =
2\mu_r \curl \vec w + \rho \vec f \\
 \nabla \cdot (\rho \vec v)&=0, \qquad \nabla \cdot \vec v =0, \\
-(c_a+c_d) \Delta \vec w
+ \rho  \jmath (\vec v \cdot \nabla )\vec w - (c_0-c_a+c_d) \nabla
(\nabla &\cdot \vec w)+ 4\mu_r \vec w =2\mu_r \curl \vec v + \rho \vec g
 \end{aligned}\end{equation}
which we assume to hold in a bounded planar domain $\Omega $, having a $C^2$
boundary $\partial \Omega$, subject to the following boundary conditions
\begin{equation}\label{e:bc_0}
 \rho = \rho _0 >0 \quad \textrm{on}~\Gamma, \qquad
\vec v=\vec v_0,~\vec w=\vec w_0 \quad \textrm{on}~\partial \Omega, \quad
\textrm{with}~ \int _{\partial \Omega } \vec v_0 \cdot \vec n\,ds =0
\end{equation}
with $\Gamma  \subset \partial \Omega$ being
 a connected arc on which $\vec v_0 \cdot \vec n < 0$, that is, $\Gamma$ is the part of $\partial \Omega $ where the fluid flows inward. 

The equations~\eqref{e:vect_f} represent conservation of linear momentum, the continuity equation, the incompressibility and the conservation of angular momentum respectively.  The unknowns are $\rho$, the denstity;
$\vec v$ and $\vec w$, the fields of velocity and rotation of particles and
the pressure, $p$. The fields $\vec f$ and $\vec g$ are, respectively,
given external sources of linear and angular momenta
densities whereas $\mu, \mu_r, c_0,c_a,c_d, \jmath$ are
positive constants characterizing the medium and also
satisfying $c_0>c_a+c_d$.

Notice that this model contains the incompressible, density
dependent Navier-Stokes system as a particular case
($\mu_r =0,\vec w\equiv 0$), and it is named the nonhomogeneous micropolar fluid
model. More details on the derivation of the model as well as the physical meaning of the several parameters above may be found in~\cite{lu_99,er_66}.

Below we devised a way to solve problem~\eqref{e:vect_f}-\eqref{e:bc_0} by combining arguments used in~\cite{lu_99,frolov93}.  This is a new result in theory and, although its proof follows fairly well-grooved lines, we believe it may serve as a first step towards a better understanding of similar flows in different geometries such as infinite or semi-infinite pipes.

Begining with Frolov's stream-function approach allows us
to drop the continuity equation and the boundary condition on $\rho$, cf.
Section~\ref{s:calc}. Then we consider an auxiliary linear problem, for the rotational velocity $\vec w$, which we solve with the Lax-Milgram lemma, see
Section~\ref{s:laxmilgram}. Next we solve a problem for the
translational velocity, $\vec v$, by defining a suitable operator and using the
Leray-Schauder principle, cf. Section~\ref{s:LS_thm}. Thus
we get a pair $\vec v, \vec w$ as a weak solution to our problem. We conclude our paper with some remarks
concerning the solution so obtained.

\section{Notations and statement of the main result}

We shall now introduce some notations and clarify what
is meant by~\eqref{e:vect_f} to hold in a planar domain.
Given $\vec u=(u_1,u_2,u_3),\vec v=(v_1,v_2,v_3)$,   we denote $[(\vec u\cdot \nabla )\vec v]_j
= \vec u \cdot \nabla v_j = \sum _k u_k \partial _{x_k} v_j$, $j=1,2,3$.
Loosely speaking, a planar flow may be sought as a ``slice" of a 3D one, that is, the flow takes place in a cross-section $x_ 3=\text{const.}$ of a 3D domain. Thus, the functions involved are assumed to be independent of the $x_3$ variable and the axes of rotation of the particles of the fluid are assumed to be perpendicular to the plane of the flow. This way, we regard
$$
\begin{gathered}
\vec v=(v_1(x_1,x_2),v_2(x_1,x_2),0),
p=p(x_1,x_2), \vec w=(0,0,w_3(x_1,x_2)), \\
\vec f=(f_1(x_1,x_2),f_2(x_1,x_2),0), \vec g=(0,0,g_3(x_1,x_2)),  
\end{gathered}$$
 and
 write $\nabla ^\perp  \psi= (-\partial
_{x_2}\psi, \partial _{x_1}\psi)$ and $\curl (\phi _1, \phi
_2) = \partial _{x_1}\phi _2 - \partial _{x_2}\phi _1$ so that
the system~\eqref{e:vect_f}
may be written componentwise as
\begin{equation}\label{e:cntr_f}
\begin{aligned}
-(\mu+\mu_r) \Delta \vec v_j + \rho \vec v \cdot \nabla v_j + \partial
_{x_j} p &=
(-1)^{j-1}2\mu_r \partial _{x_j} w_3 + \rho f_j,~j=1,2\\
 \qquad \nabla \cdot (\rho \vec v)=0, \,
\nabla \cdot \vec v &=0  \\
-(c_a+c_d) \Delta
w_3 + \rho \jmath \vec v \cdot \nabla w_3 + 4\mu_r w_3 &= 2\mu_r \curl \vec v + \rho
g_3, ~\text{in}~\Omega.
\end{aligned}
\end{equation} From now on, we shall adopt $\jmath \equiv 1,~\sigma = \mu + \mu_r$ and $\kappa = c_a+c_d$ to shorten the equations.

We use standard notations regarding Sobolev spaces modelled in $L^q(\Omega)$, $W^{k,q}(\Omega),k \ge 0, q>1,$ and their norms $\|\cdot \|_{W^{k,q}}$. The same goes to the trace spaces, $W^{k-1/q,q}(\partial \Omega ),k\ge 0, q> 1,$ and their norms $\|\cdot \|_{W^{k-1/q,q}}$. For $q=2$, we write 
$W^{k,2}(\Omega)=H^k(\Omega)$ and $W^{1/2,2}(\partial \Omega)=H^{1/2}(\partial \Omega), ~k\ge0,$ as usual.

By $\mathcal V$ we denote the set of divergence-free vector
fields $\varphi = (\varphi _1, \varphi _2)$ such that $\varphi \in C^\infty
_0 (\Omega)$, $V$ is the closure of $\mathcal V$ in the $H^1$-norm and
$\mathbf H=\{\varphi \in H^1 \mid \dvg \varphi =0~ \text{in}~
\Omega \}$. By $C^{m,\beta}(\Omega)$ we
denote the set of all $m$ times continuously differentiable functions in
$\Omega$ whose
$m$-th order derivatives are H\"older continuous with exponent $\beta \in
(0, 1)$. As above $\|\cdot \|_{C^\beta}$ stands for the $C^{0,\beta}(\Omega)$-norm.

Our main result then reads

\begin{theorem}\label{t:thm1}
Let $\vec f, \vec  g \in L^2(\Omega )$, $\rho _0 \in C^{0, \beta}(\Gamma),~ 0<\beta <
1, $ and $\vec v_0, \vec  w_0 \in H^{1/2}(\partial \Omega)$ be given
satisfying~\eqref{e:bc_0}. There exists a weak solution $\rho \in
C^{0,\alpha}(\overline \Omega), \alpha < \beta  , \vec v \in \mathbf H$, $\vec w \in
H^1(\Omega)$, of
system~\eqref{e:cntr_f} in the sense of Definition~\ref{t:defn1} below,
provided that $\mu , \kappa$ are sufficiently large so that
$
\min \bigl\{ \mu , 2\kappa \bigr\} > C\|\eta \|_{L^\infty}
\|\vec w_0\|_{H^{1/2}}.
$
\end{theorem}

\section{Weak formulation of the
problem}\label{s:calc}

For a given a divergence-free vector field $\vec v=(v_1,v_2)$ in $\Omega$
there exists $\phi : \Omega \to \mathbb R$ such that $\vec v=\nabla ^\perp
\phi$.
In addition, denoting by $\vec \tau, \vec n $ the unit tangent and outward normal
fields on $\partial \Omega$ and bearing in mind conditions~\eqref{e:bc_0}, the assumption $\nabla ^\perp \phi =
\vec v_0$ on $\Gamma$ amounts to
$$
\frac{\partial \phi}{\partial \vec n} = \vec v_0 \cdot \vec \tau, \qquad \frac{\partial
\phi}{\partial \vec  \tau} = -(\vec v_0 \cdot \vec n), \quad x \in \Gamma.
$$

Thus, the boundary values of $\phi$ may be obtained upon integration,
with respect to the arc length, from
a point $\overline{x}\in \Gamma$, $\phi (x)=- \int ^{x}_{\overline{x}} \vec v_0
\cdot \vec  n\,ds,~x \in \Gamma$. Moreover $\phi \in
H^{3/2}(\Gamma)\subset C(\Gamma),$ as $\vec v_0 \cdot \vec n
\in H^{1/2}(\Gamma)$.
From $\vec v_0 \cdot \vec n < 0$ on $\Gamma$ we see that $\phi$ is increasing on
$\Gamma$ and we may speak of its inverse $\phi ^{-1}:\phi(\Gamma)
\subset \mathbb{R} \to  \Gamma$. We may then define $\tilde{\eta} (y)=\rho
_0 (\phi ^{-1}(y)),~y \in \phi (\Gamma)\subset \mathbb R$ and
extend it to $\mathbb R$ as a strictly positive function $\eta \in
C^{0,\beta }(\mathbb R), \beta  <1 ,$ such that $\eta
(\psi (x)) =
\rho _0(x), x \in \Gamma$, whenever $\nabla
^\perp \psi = \vec v_0$ on $\Gamma $.

We fix this $\eta$ meeting the above requirements. For
sufficiently smooth $\eta, \psi$ we have  $\nabla \cdot [\eta (\psi) \nabla
^\perp
\psi ]= \eta ^\prime (\psi)\nabla \psi \cdot \nabla ^\perp \psi \equiv 0$,
in $\Omega$. A weak
version of it may be proved by regularizing $\eta$ if its assumed to be merely continuous.
From all the above facts, we assume that the density has the form $\rho = \eta
(\psi)$, where $\nabla ^\perp \psi = \vec v$, whence, the continuity
equation~\eqref{e:cntr_f}$_2$ may be dropped as well as the boundary condition~\eqref{e:bc_0}$_1$ on $\rho$, see~\cite{frolov93}.

We follow~\cite{illa01} and denote by $N: \mathbf H \to H^2(\Omega)$ the continuous operator assigning to each divergence-free vector field, $\vec u$, in $\Omega$ its stream-function $\psi = N\vec u$. Actually such a stream-function is
determined up to an arbitrary additive constant which we take to be
zero with no loss of generality, see e.g.~\cite[Lemma 2.5, Chapter
1]{tem01}
or~\cite[Theorem~4]{hey96}.

Let $\vec v_0, \vec w_0 \in H^{1/2}(\partial \Omega)$ satisfy~\eqref{e:bc_0}$_3$.
From the trace theorem follows the existence of $  b \in H^1(\Omega)$ with
\begin{equation}\label{e:bdry_w}
 b|_{\partial \Omega} = \vec w_0  \text{~and the estimate~} \|  b\|_{H^1} \le C
\|\vec w_0\|_{H^{1/2}},
\end{equation} holds for some absolute constant $C$. Also, for each $\delta >0$
fixed, there exists the
so-called Leray-Hopf extension of $\vec v_0$, that is, a vector
field $\vec a\in H^2(\Omega)$ satisfying,
\begin{equation}\label{e:lhopf_p}
\begin{gathered}
\vec a |_{\partial \Omega } = \vec  v_0, \qquad \int _{\Omega} \vec a^2 \cdot \varphi ^2
\,dx \leq
\delta ^2 \int _{\Omega}
|\nabla \varphi |^2 \,dx, \qquad \varphi \in V \\
\dvg \vec a = 0, \text{~in~} \Omega, \qquad  \vec a(x)=0, \quad d(x, \partial
\Omega) > \varepsilon, \text{~for a fixed~} \varepsilon >0.
\end{gathered}
\end{equation}
For a construction of the above cut-off functions, the reader is referred
to~\cite{Lady}. From this point on we shall omit the domain of integration
since no boundary integrals appear in the subsequent computations.

Next we introduce new unknowns $\vec{u}=\vec{v}-\vec{a}$, $\underline{ w}=w-b $, $\vec{u}\in V$ and $\underline{w}\in H^1_0(\Omega)$,
 satisfying the
following system of equations in $\Omega$
\begin{equation}\label{bvp_u}
\begin{aligned}
- \sigma\Delta \vec u +\rho [(\vec u \cdot \nabla )\vec u + (\vec a \cdot \nabla )\vec u+ (\vec u
\cdot \nabla )\vec a] +
 \nabla p &=  - 2\mu_r \nabla^\perp \underline{ w} + \rho
\vec f + \tilde F \\
 \nabla \cdot (\rho \vec u)=0, \qquad \nabla \cdot \vec u &=0,\\
-\kappa \Delta \underline{w} + \rho   [\vec u \cdot \nabla \underline{ w} +
\vec u \cdot \nabla  b + \vec a \cdot \nabla \underline{ w}]+
4\mu_r \underline{w} & =  2\mu_r \curl \vec u +
\rho \vec g + \tilde G, \\
\vec u, \underline{w} = 0,~\text{on}~\partial \Omega,
\end{aligned}
\end{equation}
with $\tilde F= \tilde F(\rho, \vec a, b):=  \sigma \Delta \vec a -\rho (\vec a \cdot
\nabla) \vec a -
2\mu_r \nabla ^\perp b$ and
$\tilde G=\tilde G(\rho, \vec a,b):= \kappa\Delta b - \rho   (\vec a \cdot \nabla) b-
4\mu_r b + 2\mu_r \curl \vec a$.

To solve the boundary value problem~\eqref{bvp_u} we look
for
$\rho= \eta (N(u+a))$ and set
$
F=F(\vec u,\vec a,b):=\tilde F(\eta (N[\vec u+\vec a]),\vec a,b)
\text{~and~}
G=G(\vec u,\vec a,b):=\tilde G(\eta (N[\vec u+\vec a],\vec a,b).
$
From the imbedding $H^2 (\Omega )\subset C^{0,\theta}(\overline{\Omega}),
\theta \in (0,1),$ and $\eta \in C^{0,\beta }(\mathbb R)$ we see
$\rho \in C^{0,\alpha}(\overline{\Omega})$, for $\alpha = \beta  \theta <
\beta $.

\begin{definition}\label{t:defn1}
We call a pair $\vec u \in V, \underline w \in H^1_0 (\Omega)$ a weak solution
of system~\eqref{bvp_u} if the integral identities
\begin{equation}\label{e:weakbvp_u}
\begin{aligned}
 \sigma\int \nabla \vec u \cdot \nabla \varphi \,dx = & \int \eta (N(\vec u+\vec a)) [\vec u
\cdot \nabla \varphi \cdot \vec u + \vec a \cdot \nabla \varphi \cdot \vec u+ \vec u
\cdot \nabla \varphi \cdot \vec a] \,dx \\ &{}-
2\mu_r \int \nabla^\perp \underline w \cdot \varphi \,dx + \int \eta (N(\vec u+\vec a))
\vec f \cdot \varphi \,dx + \int F \varphi \,dx \\
\kappa \int \nabla \underline w \nabla \xi \,dx = &  \int \eta (N(\vec u+\vec a))
[(\vec u \cdot \nabla \xi )\underline w + (\vec u \cdot \nabla \xi ) b + (\vec a \cdot
\nabla \xi ) \underline w ]\,dx \\ &{}-
2\mu_r \int (2\underline w - \curl \vec u )\xi \,dx +
\int \eta (N(\vec u+\vec a)) g \xi \,dx + \int G \xi \,dx
\end{aligned}
\end{equation}
hold for all $\varphi \in V, \xi \in H^1_0(\Omega)~\text{and}~ F, G$ as above.
\end{definition}

It is readily seen that a weak solution $\vec u, \underline{w}$ of
system~\eqref{bvp_u} in the sense of the above definition yields a weak
solution $\vec v=\vec u+\vec a, w=\underline{w}+b$ to the original
problem~\eqref{e:vect_f}-\eqref{e:bc_0}. Indeed, recalling~\eqref{e:bdry_w}
and~\eqref{e:lhopf_p}, we see that the boundary
conditions~\eqref{e:bc_0}$_{2}$ hold in the sense of traces. Moreover the
integral identities~\eqref{e:weakbvp_u} for $\vec u=\vec v-\vec a$, $\underline w =
w-b$ imply analogous ones for $\vec v,w$.
 
It is worth remarking that the recovery of the pressure is a standard
matter and it follows from the knowledge of $\rho, \vec v, w$, see
e.g.~\cite{tem01}. To be precise, we state it below as

\begin{theorem}\label{t:pres}
Let $\rho \in C^{0,\beta}(\overline{\Omega}),0< \beta< 1, \vec v \in \mathbf{H}, w\in H^1 (\Omega)$ be a weak solution of the problem~\eqref{e:vect_f}-\eqref{e:bc_0}. Then there is a $p \in L^2(\Omega)$, such that $\nabla p \in L^2(\Omega)$ and for all $\psi \in C^{\infty}_0(\Omega)$ it holds that
$$
 \sigma\int \nabla \vec v \cdot \nabla \psi \,dx - \int \rho \bigl((\vec v
\cdot \nabla )\psi \cdot \vec v - \vec f\cdot \psi \bigr)\,dx +
2\mu_r \int \nabla^\perp w \cdot \psi \,dx =\int p \dvg \psi \,dx.
$$

\end{theorem}

\section{Proof of Theorem~\ref{t:thm1}}
We split this proof in three steps.
\subsection{Auxiliary problem}\label{s:laxmilgram}
In this section we consider the following auxiliary problem:
(A) given $\vec v \in \mathbf H$, find $\underline w\in H^1_0(\Omega)$ such that
the identity
$$
-\kappa \Delta \underline{w} +
\eta(N\vec v)\vec v \cdot \nabla \underline{w} + 4\mu _r
\underline{w} = 2\mu _r \curl \vec v+ \eta(N\vec v)g + G(\vec v,b), ~\text{in}~\Omega
$$
holds in the sense of distributions, with $G(\vec v,b)=\kappa
\Delta b -
\eta(N\vec v)\vec v \cdot \nabla b - 4\mu _r b$. Existence of an unique $\underline w
\in H^1_0(\Omega)$ solving problem (A) follows from the Lax-Milgram lemma.
In fact, the problem (A) amounts to find $\underline{w} \in H^1_0(\Omega)$
such that
\begin{equation}\label{e:weakA}
B_v[\underline w, \xi]= \int \Bigl( 2\mu _r \curl \vec v + \eta (N\vec v)g +
G(\vec v,b) \Bigr) \xi \,dx ,\quad \text{for all}~\xi \in H^1_0(\Omega)
\end{equation}
and $B_v:H^1_0 (\Omega)\times H^1_0
(\Omega) \to \mathbb R$ being defined as
$$
B_v[\chi, \xi]=\int \bigl[\kappa \nabla \chi \cdot
\nabla \xi +
\eta(N\vec v)\vec v \cdot \nabla \chi \xi + 4\mu _r
\chi \xi \bigr]\,dx,\quad \text{for all}~\chi,\xi \in H^1_0(\Omega).
$$ %
According to standard
estimates $B_v$ is readily seen to be continuous.
 
It is also coercive as for all $\vec v\in \mathbf H$, $\underline w \in H^1_0(\Omega
)$ we have $\int \eta(N\vec v)\vec  v \cdot \nabla \underline w
\:\underline w \,dx=0,$ whence
$$
B_v[\underline{w},\underline{w}]=\kappa\|\nabla \underline{w}\|^2_{L^2} +
4\mu_r \| \underline{w}\|^2_{L^2} \ge
\min \{ \kappa, 4\mu_r\}\|\underline{w}\|^2_{H^1},
$$
for all $\underline w \in H^1_0(\Omega)$.

In addition, the right-hand-side of~\eqref{e:weakA} is a continuous
form in $H^1_0(\Omega)$,
$$
\left|\int \Bigl(2\mu _r \curl v + \eta (N\vec v)g +
 G(\vec v,b)\Bigr)\xi \,dx \right| \le  C \|\xi\|_{H^1},
$$ 
for all $\xi \in H^1_0(\Omega)$, with a constant $C$ depending on $c_a,c_d,
\mu _r,\|\eta\|_{L^\infty},
\|g\|_{L^2}$, $\|\vec w_0 \|_{H^{1/2}}$, $\|\vec v\|_{H^1}$ and $\Omega$.
The Lax-Milgram lemma assures the existence of an unique
$\underline{w} \in H^1_0(\Omega)$ solving problem (A).

 We point out for future reference that
 the following estimate holds for $\underline w$:
 \begin{equation}\label{e:estw}
 \begin{aligned}
 \kappa \|\nabla \underline w \|^2 _{L^2}+4\mu_r  & \| \underline w\|^2
_{L^2}
 \le  2\mu_r \|\nabla \vec v\|_{L^2}\|\underline w\|_{L^2}\\&{}+\|\eta
 \|_{L^\infty}
 \|g\|_{L^2}\|\underline
 w\|_{L^2} + C \|\eta
 \|_{L^\infty}\|\vec w_0\|_{H^{1/2}}\|v\|_{L^4}  \|\underline
 w\|_{L^2} \\&{} +C\kappa \|\vec w_0\|_{H^{1/2}}\|\nabla \underline
 w\|_{L^2} +4\mu_r C
  \|\vec w_0\|_{H^{1/2}}\|\underline
 w\|_{L^2}.
 \end{aligned}
 \end{equation}

\subsection{Problem for $u$}\label{s:LS_thm}
Our goal is to obtain $\vec u \in V$ solving~\eqref{e:problu_dif} as a fixed
point of the operator $\mathcal A
$, to be defined below. To this end we shall apply the
Leray-Schauder principle, which requires the operator
$\mathcal A
$ to be completely continuous and also that every possible
solutions of $\vec u = \lambda \mathcal A \vec u$, $\lambda  \in [0,1],$ are
uniformly
bounded, see~\cite{lu_99,Lady}.

Define $\mathcal A : V \to V$ as follows: for $\vec u \in
V$, let $\underline w \in H ^1 _0(\Omega)$ denote the solution of
problem (A) corresponding to $\vec v=\vec u+\vec a$ and
consider $\mathcal A\vec u$ given by the identity
\begin{equation}\label{e:foru}
\begin{aligned}
 \sigma\int \nabla (\mathcal A\vec u )\cdot \nabla \varphi \,dx= & \int \eta
(N[\vec u+\vec a])[(\vec u\cdot \nabla )\varphi
\cdot
\vec u + (\vec u\cdot \nabla )\varphi \cdot \vec a+(\vec a\cdot \nabla )\varphi
\cdot \vec u] \,dx \\ & {}- \int \Bigl( 2\mu_r \nabla ^\perp \underline w
- \eta (N[\vec u+\vec a])f \Bigr) \cdot \varphi \,dx + \int F\cdot \varphi \,dx,
\end{aligned}
\end{equation} for all $\varphi \in V$.
The well-definiteness of $\mathcal A$ follows from Riesz theorem. Indeed,
the rigth-hand side of~\eqref{e:foru} is a continuous form on $V$ owing to
the properties of $\vec a$, cf.~\eqref{e:lhopf_p}, and some elementary
estimates.
 
 Now we check $\mathcal A$ is completely continuous. For
$\vec u^i\in V, i=1,2$ let us define $\vec v^i = \vec u^i +\vec a \in \mathbf H$ and consider
 $\underline w ^i \in H^1 _0 (\Omega)$ the solutions of problem
(A) corresponding to $\vec v^i$. Then for $\varphi \in V,~i=1,2$,
$$
\begin{aligned}
 \sigma\int \nabla (\mathcal A\vec u^i ) \cdot \nabla \varphi \,dx = & \int \eta
(N[\vec u^i+\vec a])(\vec u^i\cdot \nabla)
\varphi \cdot \vec u^i \,dx \\
& {}+ \int \eta (N[\vec u^i+\vec a])[(\vec u^i\cdot \nabla )\varphi \cdot \vec a + (\vec a\cdot
\nabla )\varphi
\cdot \vec u^i] \,dx \\ & {}- \int \Bigl(2\mu_r \nabla ^\perp  \underline w ^i
-
\eta (N[\vec u^i+\vec a])f \Bigr)\cdot \varphi \,dx +\int F _i\cdot \varphi \,dx,
\end{aligned}
$$ with $F_i=F(\eta (N(\vec u^i+\vec a),\vec a,b)$.

Subtracting these two identities we get
\begin{equation} \label{e:problu_dif}
 \sigma\int \nabla (\mathcal A\vec u^2 - \mathcal A\vec u^1) \cdot \nabla \varphi
\,dx = I_1 + I_2 + I_3 +I_4+I_5+I_6,
\end{equation}
for
$$
\begin{aligned}
I_1 &= \int \Bigl(\eta (N[\vec u^2+\vec a])(u^2\cdot \nabla )\varphi \cdot \vec u^2 -\eta
(N[\vec u^1+\vec a])(\vec u^1\cdot \nabla )\varphi \cdot \vec u^1\Bigr)\,dx \\
I_2 &= \int \Bigl(\eta (N[\vec u^2+\vec a])(\vec u^2\cdot \nabla )\varphi \cdot \vec a - \eta
(N[\vec u^1+\vec a])(\vec u^1\cdot \nabla )\varphi \cdot \vec a\Bigr)\,dx \\
I_3 &= \int \Bigl(\eta (N[\vec u^2+\vec a])(\vec a\cdot \nabla )\varphi \cdot \vec u^2 - \eta
(N[\vec u^1+\vec a])(\vec a\cdot \nabla )\varphi \cdot \vec u^1\Bigr)\,dx \\
I_4 &=- \int  2\mu_r \nabla ^\perp (\underline w^2-\underline w^1 )\cdot
\varphi \,dx \\
I_5 &= +\int \bigl[\eta (N[\vec u^2+\vec a])-\eta (N[\vec u^1+\vec a]))\bigr]\vec f\cdot \varphi
\,dx \\
I_6 &= - \int \Bigl[\eta (N[\vec u^2+\vec a]))-\eta (N[\vec u^1+\vec a]))\Bigr] (\vec a \cdot \nabla
)\vec a \cdot \varphi \,dx.
\end{aligned}
$$

We now show it is possible to bound each of the $I_k,~k=1,\dots, 6$, by a
constant times $\|\vec u^2-\vec u^1\|_{L^4}$, proceeding as follows. First notice
that
$$
\begin{aligned}
I_1 &= \int \Bigl(\eta (N[\vec u^2+\vec a])-\eta (N[\vec u^1+\vec a])\Bigr)(\vec u^2\cdot \nabla
)\varphi \cdot \vec u^2 \,dx \\ & \ \ \ {}+\int \eta (N[\vec u^1+\vec a]) \Bigl((\vec u^2\cdot
\nabla )\varphi \cdot \vec u^2 - (\vec u^1\cdot \nabla )\varphi \cdot \vec u^1 \Bigr)\,dx
\\
&= \int \Bigl(\eta (N[\vec u^2+\vec a])-\eta (N[\vec u^1+\vec a])\Bigr)(\vec u^2\cdot \nabla
)\varphi \cdot \vec u^2 \,dx \\ & \ \ \ {}+\int \eta (N[\vec u^1+\vec a]) \Bigl\{[(\vec u^2 -
\vec u^1) \cdot \nabla ]\varphi \cdot \vec u^2 + (\vec u^1\cdot \nabla )\varphi \cdot
(\vec u^2-\vec u^1) \Bigr\}\,dx=I_{11}+I_{12}.
\end{aligned}
$$
H\"older, Young and imbedding inequalities imply
$$
\begin{aligned}
|I_{11}| &\le \|\eta\|_{C^\alpha} \int |N(\vec u^2-\vec u^1)|^\alpha \Bigl|\sum
_{j,k}  u^2 _k \partial _{x_k}\varphi _j  u^2 _j \Bigr|\,dx \\
& \le \|\eta\|_{C^\alpha} \sum _{j,k} \Bigl(\int (\partial _{x_k}\varphi
_j)^2  \,dx\Bigr)^{1/2} \Bigl(\sum _{j,k}\int |N(\vec u^2-\vec u^1)|^{2\alpha}
[ u^2 _k  u^2 _j]^2 \,dx\Bigr)^{1/2}  \\
& \le C\|\eta \|_{C^\alpha }\|\nabla \varphi \|_{L^2}\|N(\vec u^2 -
\vec u^1)\|^\alpha _{W^{1,4\alpha}}\|\vec u^2\|^2_{L^8} \\ &=C\|\eta \|_{C^\alpha
}\|\nabla \varphi \|_{L^2}\|\vec u^2 - \vec u^1\|^\alpha _{L^{4\alpha}}
\|\vec u^2\|^2_{L^8}, \\
|I_{12}| & \le \|\eta\|_{L^\infty}\|\nabla \varphi
\|_{L^2}(\|\vec u^2\|_{L^4}+\|\vec u^1\|_{L^4})\|\vec u^2 -\vec  u^1\| _{L^4}.
\end{aligned}
$$

Hence
\begin{equation}\label{e:u21l4_i1}
\begin{aligned}
|I_1| \le C\|\eta \|_{C^\alpha
}\|\nabla \varphi \|_{L^2}\|\vec u^2 & - \vec u^1\|^\alpha _{L^{4\alpha}}
\|\vec u^2\|^2_{L^8}\\
& {}+ \|\eta\|_{L^\infty}\|\nabla \varphi
\|_{L^2}(\|\vec u^2\|_{L^4}+\|\vec u^1\|_{L^4})\|\vec u^2 - \vec u^1\| _{L^4}.
\end{aligned}
\end{equation}
Likewise we get the bounds
\begin{equation}\label{e:u21l4_i2}
\begin{gathered}
\begin{aligned}
|I_2| \le \|\eta\|_{L^\infty}\|\nabla  \varphi \|_{L^2}\|\vec u^2- & \vec u^1
\|_{L^4}\|\vec a\|_{L^4} \\ &{}+ C\|\eta \|_{C^\alpha
}\|\vec u^2\|_{L^8}\|\vec a\|_{L^8}\|\nabla \varphi \|_{L^2}\|\vec u^2 - \vec u^1\|^\alpha
_{L^{4\alpha}}
\end{aligned} \\
\begin{aligned}
|I_3| \le \|\eta\|_{L^\infty}\|\nabla \varphi \|_{L^2}\|\vec u^2- & \vec u^1
\|_{L^4}\|\vec a\|_{L^4} \\ &{}+ C\|\eta \|_{C^\alpha
}\|\vec u^1\|_{L^8}\|\vec a\|_{L^8}\|\nabla \varphi \|_{L^2}\|\vec u^2 - \vec u^1\|^\alpha
_{L^{4\alpha}}
\end{aligned} \\
|I_5| \le  \|\eta \|_{C^ \alpha}\|\vec u^2-\vec u^1\|^\alpha _{L^{4\alpha
}}\|\vec f\|_{L^2}\|\varphi \|_{L^4} \\
\begin{aligned}
|I_6| \le  C\|\eta
\|_{C^\alpha } \|\vec a\|^2_{L^8}\|\nabla \varphi \|_{L^2}\|\vec u^2 - \vec u^1\|^\alpha
_{L^{4\alpha}}.
\end{aligned}
\end{gathered}
\end{equation}

The term $I_4$ requires estimating $\|\underline w^2-\underline
w^1\|_{H^1}$ in terms of $\|\vec u^2-\vec u^1\|_{L^4}$, a task we now perform.
As $\underline w ^i,~i=1,2$ solve problem (A) we find that the following
identity
$$
\int \kappa \nabla (\underline w ^2 - \underline w ^1)\cdot  \nabla \psi +
4\mu_r (\underline w ^2 - \underline w ^1) \psi \, dx = J_1 +J_2 +J_3,
$$
holds for all $\psi \in H^1_0(\Omega)$, where
$$
\begin{aligned}
J_1 &= \int \Bigl(\eta (N[\vec u^2+\vec a])(\vec u^2+\vec a) \cdot \nabla \underline w^2 -\eta
(N[\vec u^1+\vec a])(\vec u^1+\vec a)\cdot \nabla \underline w^1 \Bigr) \psi \,dx \\
J_2 &= \int \Bigl(2\mu_r \nabla \times (\vec u^2 - \vec u^1) + \bigl[ \eta
(N[\vec u^2+\vec a])(\vec u^2+\vec a) - \eta (N[\vec u^1+\vec a])(\vec u^1+\vec a)\bigr]g \Bigr) \psi \,dx \\
J_3 &= \int \Bigl(\eta (N[\vec u^2+\vec a])\vec u^2 -\eta (N[\vec u^1+\vec a])\vec u^1 \Bigr)\cdot \nabla
b \psi \,dx.
\end{aligned}
$$

By arguing as in the estimations~\eqref{e:u21l4_i1}-\eqref{e:u21l4_i2} we
find
\begin{align*}
|J_1| & \le C \Bigl( \|\eta \|_{C^\alpha} \| \vec u^2 - \vec u^1
\|_{L^{4\alpha}}^{\alpha}\|\vec u^1+\vec a\|_{L^8} +
\|\eta\|_{L^\infty}\|\vec u^2-\vec u^1\|_{L^4} \Bigr) \|\nabla \underline w^2\|_{L^2}
\|\psi \|_{H^1} \\
|J_2| & \le C \Bigl( 2\mu_r\text{meas}(\Omega)^{1/4}\|\vec u^2-\vec u^1\|_{L^{4}}+
\|\eta\|_{C^{\alpha}}\|g\|_{L^{2}}\|\vec u^2-\vec u^1\|_{L^{4\alpha}}^\alpha \Bigr)
\|\psi \|_{H^1} \\
|J_3| &  \le C \Bigl(
\|\eta\|_{C^{\alpha}}\|\vec u^2-\vec u^1\|_{L^{4\alpha}}^{\alpha}\|\vec u^1+\vec a\|_{L^8} +
\|\eta\|_{L^{\infty }}\|\vec u^2-\vec u^1\|_{L^{4}} \Bigr) \|\nabla b\|_{L^2} \|\psi
\|_{H^1}.
\end{align*}
Taking $\psi = \underline w^2 -\underline w^1 $ and invoking the
boundedness of $\Omega$, we get
$$
\begin{aligned}
\|\underline w^2 - \underline w^1 \|_{H^1} &\le C \Bigl(\|\nabla \underline
w^2\|_{L^2}+ \|\nabla b\|_{L^2}  \Bigr)
\Bigl\{\|\eta\|_{C^{\alpha}}\|\vec u^1+\vec a\|_{L^8}\|\vec u^2-\vec u^1\|_{L^{4}}^{\alpha} \\
& {}+\|\eta\|_{L^\infty}
\|\vec u^2-\vec u^1\|_{L^{4}}\Bigr\} + (2\mu_r\text{meas}(\Omega)^{1/4}+
\|\eta\|_{C^\alpha }\|g\|_{L^2})\|\vec u^2-\vec u^1\|_{L^{4}},
\end{aligned}
$$ for some constant $C>0$ depending on $\kappa, \mu_r, \Omega$. At last,
bearing in mind the triangle inequality, we find
\begin{equation}\label{e:est_i4}
\begin{aligned}
|I_4| &\le 2\mu_r C \Bigl(\|\nabla \underline w^2\|_{L^2}+ \|\nabla
b\|_{L^2}  \Bigr)\|\eta\|_{C^{\alpha}}(\|\vec u^1\|_{L^8}+\|\vec a\|_{L^8}
)\|\vec u^2-\vec u^1\|_{L^{4}}^{\alpha}\|\varphi \|_{H^1} \\
&{}+2\mu_r
C(2\mu_r\text{meas}(\Omega)^{1/4}+\|\eta\|_{C^\alpha}\|g\|_{L^2})
\|\vec u^2-\vec u^1\|_{L^{4}}\|\varphi \|_{H^1}.
\end{aligned}
\end{equation}

Collecting inequalities~\eqref{e:u21l4_i1}-\eqref{e:est_i4}  for $\varphi =
\mathcal A \vec u^2 - \mathcal A \vec u^1$, together with  $\|\vec u^2-\vec u^1\|_{L^4} \leq
1$ and in view of~\eqref{e:bdry_w} and~\eqref{e:lhopf_p}, we conclude
\begin{equation*}
\begin{aligned}
 \sigma \| \nabla (\mathcal A\vec u^2 - \mathcal A\vec u^1) \|_{H^1}  \le \|\vec u^2 -&
\vec u^1\|^\alpha _{L^{4}}  \times \left\{ C\|\eta
\|_{C^\alpha } \Bigl[\|\vec u^2\|_{L^8}
 (\|\vec u^2\|_{L^8}+C \| \vec v_0 \|_{H^{1/2}}) \right. \\
& {}+ C \| \vec v_0 \|_{H^{1/2}}(\|\vec u^1\|_{L^8}+C \| \vec v_0 \|_{H^{1/2}}) +
\|\vec f\|_{L^2}+ \|g\|_{L^{2}} \\
& {}+ (\|\vec u^1 \|_{L^8} + C \| \vec v_0 \|_{H^{1/2}} )(\|\nabla \underline
w^2\|_{L^2}
+C \| \vec w_0 \|_{H^{1/2}}) \Bigr] \\
& {}+ \|\eta\|_{L^\infty} \Bigl[\|\vec u^2\|_{L^4}+\|\vec u^1\|_{L^4} 
 +  2C \| \vec v_0 \|_{H^{1/2}}  \\
& {}+ (\|\vec u^1 \|_{L^8} + C \| \vec v_0 \|_{H^{1/2}} )(\|\nabla \underline
w^2\|_{L^2} +C \| \vec w_0 \|_{H^{1/2}}) \\
& \left. {}+ 2\mu_r\text{meas}(\Omega)^{1/4} \Bigr] \right\}.
\end{aligned}
\end{equation*}
Thus, from the compactness of the
imbedding $H^1 \hookrightarrow L^4$ and the inequality above, a weakly convergent sequence in $V$ is mapped by $\mathcal A$ into a strongly convergent sequence in $L^4$. It remains to show that all possible
solutions of the equation $\vec u = \lambda \mathcal A \vec u, \lambda \in [0,1],$
are uniformly bounded. As $\lambda = 0 $ implies $\vec u=0$, we
suppose $\lambda > 0 $ and replace $\mathcal A \vec u = \vec u/\lambda $,
$\varphi = \vec u$ in the equation~\eqref{e:foru}. Performing estimations
similar to the previous we find
\begin{equation}\label{e:unif_u}
\frac \sigma \lambda \|\nabla \vec u \|_{L^2}^2 \le \| \eta
\|_{L^\infty} \Bigl|\int (\vec u  \cdot \nabla )\vec  u \cdot \vec a\,dx\Bigr| +
2\mu_r \bigl(\|\underline w\|_{L^2}+
\|b\|_{L^2}\bigr)\|\nabla \vec u\|_{L^2} + \|\eta \|_{L^\infty
}\|\vec f\|_{L^2}\|\vec u\|_{L^2}.
\end{equation}

Notice that the term $2\mu_r \|\underline w\|_{L^2}\|\nabla \vec u\|_{L^2}$
itself is another quadratic term in $\vec u$ which we handle as follows. In view
of ~\eqref{e:estw}, with $\vec v=\vec u+\vec a$, we find
\begin{equation}\label{e:estw_u}
\begin{aligned}
 \kappa \|\nabla \underline w \|^2 _{L^2}+4\mu_r\| \underline w\|^2 _{L^2}
 & \le  2\mu_r \|\nabla (\vec u+\vec a)\|_{L^2}\|\underline w\|_{L^2}+\|\eta
 \|_{L^\infty}
 \|g\|_{L^2}\|\underline
 w\|_{L^2} \\&{}+ C \|\eta
 \|_{L^\infty}\|\vec w_0\|_{H^{1/2}}\|\vec u+\vec a\|_{L^4}  \|\underline
 w\|_{L^4}  +C\kappa \|\vec w_0\|_{H^{1/2}}\|\nabla \underline
 w\|_{L^2}\\&{}+4\mu_r
  \|\vec w_0\|_{H^{1/2}}\|\underline
 w\|_{L^2}.
 \end{aligned}
 \end{equation}

Aided by Young and H\"older inequalities we get, by summing
up~\eqref{e:unif_u} and~\eqref{e:estw_u},
$$
\begin{aligned}
\frac{\mu + \mu_r}{\lambda} \|\nabla \vec u\|^2 _{L^ 2} + \kappa \|\nabla
\underline w\|^2 _{L^ 2} & \le \|\eta \|_{L^ \infty}\Bigl|\int (\vec u \cdot
\nabla ) \vec u \cdot \vec a\,dx\Bigr| \\ &{}+\mu_r \|\nabla \vec u \|^2_{L^ 2} + C\|\eta
\|_{L^\infty} \|\vec u\|_{L^4} \|\vec w_0\|_{H^{1/2}}\|\underline w\|_{L^4}
\\&{}+2\mu_r\|b\|_{L^2} \|\nabla \vec u\|_{L^2}+ \|\eta
\|_{L^\infty}\|\vec f\|_{L^2}\|\vec u\|_{L^2} \\ &{}+2\mu_r\|\nabla
\vec a\|_{L^2}\|\underline w\|_{L^2} + \|\eta
\|_{L^\infty}\|g\|_{L^2}\|\underline w\|_{L^2} \\ &{}+C\|\eta \|_{L^\infty}
\|\vec a\|_{L^4} \|\vec w_0\|_{H^{1/2}}\|\underline w\|_{L^4} + C\kappa \|\nabla
\underline w\|_{L^2}\|w_0\|_{H^{1/2}} \\
&{}+4\mu_r\|\vec w_0\|_{H^{1/2}}\|\underline w \|_{L^2}.
\end{aligned}$$
Next, requiring $\delta >0 $ in~\eqref{e:lhopf_p}, to be such that $\delta
\|\eta\|_{L^\infty} < \mu /2$ and
estimating
$$
\|\eta \|_{L^\infty} \|\vec u\|_{L^4} \|\vec w_0\|_{H^{1/2}}\|\underline w\|_{L^4}
\le \frac{C}{2}\|\eta \|_{L^\infty} \|\vec w_0\|_{H^{1/2}} \Bigl( \|\nabla
\vec u\|^2_{L^2} +\|\nabla \underline w\|^2 _{L^2}\Bigr)
$$
we ultimately get, using~\eqref{e:bdry_w} and~\eqref{e:lhopf_p},
$$
\begin{aligned}
\Bigl(\frac{\mu}{2\lambda} - \frac{C}{2}&\|\eta \|_{L^\infty}\|\vec w_0\|_{H^
{1/2}}\Bigr) \|\nabla \vec u\|^2 _{L^ 2} + \Bigl(\kappa - \frac{C}{2}\|\eta
\|_{L^\infty}\|\vec w_0\|_{H^ {1/2}}\Bigr)\|\nabla \underline w\|^2 _{L^ 2}  \\
& \le 2\mu_r C\|\vec w_0\|_{H^{1/2}} \|\nabla \vec u\|_{L^2}+ \|\eta
\|_{L^\infty}\|\vec f\|_{L^2}\|\vec u\|_{L^2} +2\mu_r C \|
\vec v_0\|_{H^{1/2}}\|\underline w\|_{L^2} \\ &{}+ \|\eta
\|_{L^\infty}\|g\|_{L^2}\|\underline w\|_{L^2} +C\|\eta \|_{L^\infty}
\|\vec v_0\|_{H^{1/2}} \|\vec w_0\|_{H^{1/2}}\|\underline w\|_{L^4} \\ &{}+ C\kappa
\|\nabla \underline w\|_{L^2}\|\vec w_0\|_{H^{1/2}}
+4\mu_r C\|\vec w_0\|_{H^{1/2}}\|\underline w \|_{L^2}.
\end{aligned}$$
Therefore, demanding $\mu $ and $\kappa $ to be large enough so that
$$
\min \bigl\{ \mu , 2\kappa \bigr\} > C\|\eta \|_{L^\infty}
\|\vec w_0\|_{H^{1/2}},
$$
we may conclude the uniform boundedness on the norms of all
possible solutions of $\vec u=\lambda \mathcal A \vec u,~\lambda \in [0,1],$ with
respect to the parameter $\lambda$.
From this and previous steps we conclude that the Leray-Schauder principle
applies and that
problem~\eqref{bvp_u} has a weak solution. \qed

\section{Concluding remarks}

We have tacitly assumed $\partial \Omega$ to consist of a single component, that is, $\Omega \subset \mathbb R ^2$ to be a simply connected open set. A clue on how the arguments should be modified to cope with a non-simply connected $\Omega $ may be found in~\cite{frolov93}.
 
The presence of the quadratic term brought up by the equation~\eqref{e:vect_f}$_4$, for the rotational field $w$, required us to demand the viscosities $\mu, \kappa$ to be sufficiently large compared to data  and this somewhat contrasts with previous results by Frolov~\cite{frolov93} and Santos~\cite{msantos}.

As shown in~\cite[cf. Theorems~2 and 3]{frolov93},
also~\cite[Chapter~2]{lu_99}, the regularity of
$\rho , \vec v, w, p$ may be improved by increasing those of $\partial \Omega,
\rho _0, \vec v_0, \vec w_0, \vec f$ and $g$. An usual bootstrap argument should combine
regularity results for the Stokes problem and Ne\v cas' results on strongly
elliptic systems of second order~\cite{necas65}.

We also notice that the uniqueness of the above solution deserves further
investigation. Indeed we could not benefit from neither the steady-state
continuity equation~\eqref{e:vect_f}$_2$ nor the particular form of the
density to derive the required estimates.


\end{document}